\input amstex
\documentstyle{amsppt}
\loadeusm
\input boxedeps.tex
\SetRokickiEPSFSpecial
\HideDisplacementBoxes

\define\bandt{1}
\define\bandtg{2}
\define\barn{3}
\define\cm{4}
\define\drobot{5}
\define\edgar{6}
\define\edgaripfm{7}
\define\edgarmauldin{8}
\define\falconer{9}
\define\gantmacher{10}
\define\golds{11}
\define\mauldin{12}
\define\mw{13}
\define\mcmullen{14}
\define\minc{15}
\define\patzschke{16}
\define\schief{17}
\define\seneta{18}
\define\strichartz{19}
\define\tsuji{20}

\define\scr#1{{\fam\eusmfam\relax#1}}
\define\rest{{\restriction}}         
\define\rmin{r_{\roman{min}}}

\define\dmin{d_{\roman{min}}}
\define\dmax{d_{\roman{max}}}
\define\lambdamin{\lambda_{\roman{min}}}
\define\lambdamax{\lambda_{\roman{max}}}
\define\lambdamaxp{\lambda'_{\roman{max}}}
\define\lambdaminp{\lambda'_{\roman{min}}}
\define\?#1{#1}
\define\diam{\operatorname{diam\,}}
\define\dimB{\overline{\dim}_{\roman B}}
\define\dimH{\dim_{\roman H}}
\define\dist{\operatorname{dist}}
\define\SET#1#2{ \left\{\, {#1} : {#2} \,\right\} }     
\define\takes{\:}                    
\define\emptystring{\Lambda}
\define\N{\Bbb N}                    
\define\scrH{{\scr H}}     
\define\Def#1{{\bf #1\/}}   
\define\boldx{\bold x}

\catcode`\@=11

\font\tenlab=cmssbx10
\font\sevenlab=cmssbx10 at 7pt
\font\fivelab=cmssbx10 at 5pt 
\newfam\labfam
\textfont\labfam=\tenlab
\scriptfont\labfam=\sevenlab
\scriptscriptfont\labfam=\fivelab
\define\lab{\relaxnext@\ifmmode\let\next\lab@\else
 \define\next{\Err@{Use \string\lab\space only in math mode}}\fi\next}
\define\lab@#1{{\lab@@{#1}}}
\define\lab@@#1{\noaccents@\fam\labfam#1}
\catcode`\@=\active

\newcount\sectno                    
\newcount\itemno                    
\sectno=0
\define\newsection{%
  \global\advance\sectno by 1%
  \global\itemno = 0}
\define\theitem{\the\sectno.\the\itemno}
\define\Itemx{%
  \global\advance\itemno by 1%
  \noindent{\bf(\theitem)\ }}
\define\Itemd{%
  \global\advance\itemno by 1%
  \tag\theitem
  }

\topmatter
\title
A fractal dimension estimate for a graph-directed
IFS of non-similarities
\endtitle
\rightheadtext{A fractal dimension estimate}

\author G. A. Edgar
and
Jeffrey Golds
\endauthor

\affil The Ohio State University
and
Silicon Graphics, Inc.
\endaffil

\address Department of Mathematics, The Ohio State University,
Columbus, OH 43210\endaddress
\email edgar\@math.ohio-state.edu\endemail

\address \null \endaddress
\email golds\@csd.sgi.com \endemail






\abstract
Suppose a graph-directed iterated function system consists of
maps $f_e$ with upper estimates of the form
$d\big(f_e(x),f_e(y)\big) \le r_e d(x,y)$.  Then the fractal dimension
of the attractor $K_v$ of the IFS is bounded above by the dimension associated
to the Mauldin--Williams graph with ratios $r_e$.
Suppose the maps $f_e$ also have lower estimates of the form
$d\big(f_e(x),f_e(y)\big) \ge r'_e \,d(x,y)$ and that the IFS also
satisfies the strong open set condition.  Then the fractal dimension
of the attractor $K_v$ of the IFS is bounded
below by the dimension associated
to the Mauldin--Williams graph with ratios $r'_e$.
When $r_e = r'_e$, then the maps are similarities and this reduces to
the dimension computation of Mauldin \& Williams 
for that case.
\endabstract

\endtopmatter

\document

\head 0. Introduction \endhead

Fractal sets may be constructed in many different ways.  Barnsley
\cite{\barn} singled out the ``iterated function system'' method:
The fractal set $K$ is made up of parts, each of which is a
shrunken copy of the whole set.
Mauldin \& Williams \cite{\mw} provided a more general
setting, where several sets $K_v$ are involved, each of them is
made up of parts, and each part is a shrunken copy of
one of the parts (the same one or a different one).  The combinatorics
of the way in which the parts fit together is described
by a directed multigraph.  This is described in detail
below (Definition 1.5).  In addition to Mauldin and Williams,
compare ``recurrent iterated function system'' \cite{\barn\rm, Ch.~X},
``Markov self-similar sets'' \cite{\tsuji}, and ``mixed self-similar systems''
\cite{\bandt}; see also \cite{\drobot}, \cite{\strichartz}.

In the text \cite{\edgar\rm, \S 6.4} there is an exposition of the Mauldin--Williams
computation of the dimension of the attractors for a graph-directed
iterated function system consisting of similarities.  ``Similarities''
are functions $f \takes S \to T$ between metric spaces that satisfy equations
of the form
$$
	d\big(f(x),f(y)\big) = r\,d(x,y)\qquad\text{for all $x,y \in S$.}
\Itemd
$$
Here we have used $d$ for the (possibly different)
metrics in the two metric spaces.  This dimension computation involves
the spectral radius of some nonnegative matrices.  Next in the text
are two ``Exercises'' \cite{\edgar\rm, (6.4.9) and (6.4.10)}
asking for appropriate generalizations when the equalities (\?{0.1})
are replaced by inequalities of the form
$$
	d\big(f(x),f(y)\big) \le r\,d(x,y)\qquad\text{for all $x,y \in S$}
$$
or of the form
$$
	d\big(f(x),f(y)\big) \ge r\,d(x,y)\qquad\text{for all $x,y \in S$.}
$$
In this paper we provide some possible solutions for those
exercises.  The non-graph
case (or, in our language, a graph with only one node and several loops)
may be found in \cite{\falconer\rm, Theorems 9.6,~9.7}.

We also take this opportunity to make a generalization beyond
the setting used in \cite{\edgar\rm, \S 6.4}: we allow complete metric
spaces other than Euclidean space.   (Even when Euclidean space
is the main interest, it is useful to apply the results in 
ultrametric ``model'' spaces consisting of strings or sequences.)
This means that for the lower bound estimates we must use the ``strong
open set condition'' rather than the original ``open set condition''
(this has been discovered since the publication of~\cite{\edgar}:
for example \cite{\bandtg}, \cite{\patzschke}, \cite{\schief}).

The paper is arranged as follows.  Section~1 reviews the notation
for graph-directed iterated function systems.  Section~2 reviews
the Perron--Frobenius theory to be used.  Section~3 contains
the proofs.  Section~4 describes an example.

This paper is based in part on a portion of the dissertation
\cite{\golds} of the second author, written at The Ohio
State University under the direction of the first author.

\head 1. The Setting \endhead\newsection

We will describe here: the ``fractal dimensions'' that will be used,
the ``graph-directed iterated function systems''
and the fractal sets $K_u$ that they define, and the
``string models'' $E_u^{(\omega)}$ that will
be used in the investigation.  The notation will
mostly follow \cite{\edgar}, especially Section~4.3 on graph
self-similarity.

\subhead Fractal Dimensions\endsubhead
The first dimension we consider is the box-counting dimension.
See \cite{\falconer} for a more complete discussion.

Let $S$ be a metric space, and let $K\subseteq S$ be a totally
bounded set.  For each $\delta > 0$, define $N_{\delta}(K)$ to be
the smallest number of sets of diameter $\le \delta$ that can cover $E$.
Since we have postulated that $K$ is totally bounded, $N_\delta(K)$
is a finite natural number for every $\delta > 0$.  Unless $K$ is
a finite set, $N_\delta(K)$ increases to $\infty$ as $\delta$ decreases
to $0$; the rate at which $N_\delta(K)$ increases will tell us something
about the size of the set $K$.

\definition{\Itemx Definition}
The \Def{upper box-counting dimension} or \Def{Bouligand dimension} of $K$ is
$$
	\dimB K =
	\limsup_{\delta \to 0} \frac{\log N_\delta(K)}{-\log \delta}.
$$
The \Def{lower box-counting dimension} of $K$ is
$$
	\underline{\dim}_{\roman B} K =
	\liminf_{\delta \to 0} \frac{\log N_\delta(K)}{-\log \delta}.
$$
\enddefinition

Next we will review the definition of the Hausdorff dimension.
See for example \cite{\edgar\rm, p.~149}, \cite{\falconer\rm, p.~28}.

\definition{\Itemx Definition}
Let $S$ be a metric space,
let $K\subseteq S$, and let 
$\delta > 0$ be a positive number.
A collection $\scr C$ of subsets of $S$ is
\Def{$\delta$-fine} iff $\diam A \le \delta$ for all $A \in \scr C$.
The collection $\scr C$ is a \Def{cover} of $K$ iff
$K \subseteq \bigcup_{A \in \scr C} A$.
\enddefinition

\definition{\Itemx Definition}
Let $S$ be a metric space,
let $K\subseteq S$, and let $s > 0$.  For $\delta > 0$, define
$$
	\scrH^s_\delta(K) = \inf \sum_{A \in \scr C} (\diam A)^s ,
$$
where the infimum is over all countable $\delta$-fine
covers $\scr C$ of $K$.  Define
$$
	\scrH^s(K) = \sup_{\delta > 0} \scrH^s_\delta(K)=
	\lim_{\delta \to 0} \scrH^s_\delta(K) .
$$
The \Def{Hausdorff dimension} of the set $K$ is
$$
	\dimH K = \sup\SET{s}{\scrH^s(K) = \infty}
	= \inf\SET{s}{\scrH^s(K) = 0} .
$$
\enddefinition

Another fractal dimension (not defined
here) is the packing dimension $\dim_{\roman P}$.
See \cite{\edgaripfm\rm, \S 1.2} for a discussion.  In \cite{\edgar}
the packing measure is mentioned only as an afterthought, but
in \cite{\edgaripfm} its role is nearly
as important as the Hausdorff measure.

Many other fractal dimensions may be found in the literature.

It is known (\cite{\edgaripfm\rm, Corollary 1.3.5}
and \cite{\falconer\rm, Lemma 3.7})
that all of the dimensions mentioned here (and many other fractal
dimensions, as well)
are between the Hausdorff dimension and the upper box dimension:
$$
	\dimH K \le \dim K \le \dimB K,
$$
where ``$\dim$'' is any of
$\dimH, \dim_{\roman P},
\underline{\dim}_{\roman B},
\overline{\dim}_{\roman B}$.

\subhead Graphs and Iterated Function Systems\endsubhead
First, a directed multigraph $(V,E)$ should be
fixed.  The elements $v \in V$ are the
\Def{vertices} or \Def{nodes}
of the graph; the elements $e \in E$
are the \Def{edges} of the graph.  For $u,v \in V$,
there is a subset $E_{uv}$ of $E$, known as the
edges from $u$ to $v$.  Each edge belongs to exactly
one of these subsets.  We will sometimes write
$E_u = \bigcup_v E_{uv}$, the set of all edges
leaving the vertex $u$.

We will often think of the
set $E$ as a set of ``letters'' that label the
edges of the graph, so we will talk about
``words'' or ``strings'' made up of these letters.
A \Def{path} in the graph is a finite string
$\alpha = e_1 e_2 \cdots e_k$ of edges, such that
the terminal vertex of each edge $e_i$ is the initial
vertex of the next edge $e_{i+1}$.
The initial vertex of $\alpha$ is the initial vertex
of the first letter $e_1$ and the terminal vertex
of $\alpha$ is the terminal vertex of the last letter
$e_k$.  We write $E_{uv}^{(k)}$ for the set of all paths
of length $k$ that begin at $u$ and end at $v$;
and $E_u^{(k)}$ for the set of all paths of length
$k$ that begin at $u$; and
$E_u^{(*)}$ for the set of all finite paths
of any length that begin at $u$; and
$E^{(*)}$ for the set of all finite paths.
By convention, $E^{(0)}_u$ consists of a single
``empty path'' $\emptystring_u$ of length zero from node $u$ to itself.

If $\alpha \in E_{uv}^{(k)}$ and $\beta \in E_{vw}^{(n)}$,
then we write $\alpha \beta$ for the path made by concatenation
of the two strings, so that $\alpha \beta \in E_{uw}^{(k+n)}$.

A partial order may be defined on $E^{(*)}$ as follows:
write $\alpha \le \beta$ iff $\alpha$ is a prefix of
$\beta$: that is, $\beta = \alpha\gamma$ for some path
$\gamma$.  With this ordering, each  set $E_u^{(*)}$
becomes a tree.  The root of the tree $E_u^{(*)}$
is the empty path $\emptystring_u$.
The entire set $E^{(*)}$ is
a disjoint union of trees (a forest).
Two strings are \Def{incomparable} iff neither is a prefix of the other.

A nonempty path that begins and ends at the same node is called
a \Def{cycle}.

We will assume that the
graph $(V,E)$ is \Def{strongly connected}, that is, there is
a path from any vertex to any other, along the edges
of the graph (taken in the proper directions).
We will also assume that there are at least two edges leaving
each node.  (In \cite{\edgarmauldin} may be found a method to
eliminate this assumption.)

\definition{\Itemx Definition}
An \Def{iterated function system} (or IFS) directed by the graph $(V,E)$
consists of: one complete metric space $S_v$ for each vertex $v \in V$
and one function $f_e \takes S_v \to S_u$ for each edge $e \in E_{uv}$.
\enddefinition

To avoid trivialities, we will normally
assume that each metric space $S_v$ has at least two points.
For this sort of graph-directed arrangement see \cite{\mw}, \cite{\edgar},
and references in both.  Note that in some places, such as \cite{\strichartz},
the direction of the edges of the graph are the reverse of
the convention used in this paper; this may lead to consideration
of the transpose of our matrix, but of course that does not affect
the spectral radius.

\definition{\Itemx Definition}
An \Def{invariant set list} for the iterated function system $(f_e)$
consists of one nonempty set $K_v \subseteq S_v$ for each
$v \in V$ such that
$$
	K_u = \bigcup_{v \in V}\bigcup_{e \in E_{uv}} f_e\big[K_v\big]
$$
for all $u \in V$.
\enddefinition

If the functions $f_e$ satisfy appropriate conditions, then
it may be proved that there is a unique invariant set list
of nonempty compact sets.  For example \cite{\falconer\rm, Theorem~9.1},
\cite{\edgar\rm, Theorem~4.3.5}.

\subhead The Models\endsubhead
We will use some ``string models'' in our
investigation of these invariant sets.  Write
$E_u^{(\omega)}$ for the set of all {\it infinite\/}
strings, using symbols from $E$, where the initial vertex
of the first edge is $u$ and the terminal
vertex of each edge is the initial vertex of the next
edge.  These sets are naturally compact metric
spaces.  For each finite string $\alpha \in E^{(*)}$, the \Def{cylinder}
$[\alpha]$ is the set of all infinite strings
$\sigma \in E^{(\omega)}$ that begin with $\alpha$.
Then the set $\big\{\,[\alpha] : \alpha \in E_u^{(*)}\,\big\}$
is the base for the topology on $E_u^{(\omega)}$.
For $\sigma \in E_u^{(\omega)}$ and a positive
integer $k$, the \Def{restriction}
$\sigma \rest k$ is the finite string made up of the
first $k$ letters of $\sigma$.
The same notation $\alpha \rest k$ is used for
finite strings $\alpha$ when $k$ is less than the
length of $\alpha$.  As a special case of this:
if $\alpha$ has length $k$, then the \Def{parent}
of $\alpha$ is $\alpha^{-} = \alpha \rest (k-1)$,
obtained by omitting the last letter of $\alpha$.

Suppose $(E,V)$ directs an iterated function system $(f_e)$
consisting of contractions.
Then we may define an ``action'' of $V$ starting from the
maps $f_e$.  If $\alpha = e_1 e_2 \cdots e_k \in E_{uv}^{(k)}$,
then the composition
$f_\alpha = f_{e_1}\circ f_{e_2}\circ\cdots\circ f_{e_k}$
is a function from $S_v$ to $S_u$.  We will usually
write simply $\alpha(x)$ for $f_\alpha(x)$.

There is a
\Def{model map} $h_u \takes E_u^{(\omega)} \to S_u$
for each $u$, defined so that
$h_u(\sigma)$ is the unique element of the set
$$
	\bigcap_{k = 1}^\infty (\sigma\rest k)K_{v_k}.
$$
(Of course, here $v_k$ is the terminal vertex
for the $k$th letter of $\sigma$.)
Then clearly $K_u = h_u\big[E_u^{(\omega)}\big]$,
since the sets $h_u\big[E_u^{(\omega)}\big]$ satisfy the
defining conditions for an invariant set list of
the IFS $(f_e)$.
If $h_u(\sigma) = x$, then we say that the string
$\sigma$ is an \Def{address} of the point
$x$.

\subhead Mauldin--Williams Graphs\endsubhead
The computations to be done here can be stated in terms of
a combinatorial structure that assigns a positive number to
each edge of the graph.

\definition{\Itemx Definition}
A \Def{Mauldin--Williams Graph} is a directed multigraph
$(V, E)$ together with a positive number $r_e$ for each
edge $e \in E$.
\enddefinition

If $\alpha = e_1 e_2 \cdots e_k \in E^{(k)}$ is a path, we
define $r(\alpha) = \prod_{i=1}^k r_{e_i}$; thus
$r(\alpha\beta) = r(\alpha)r(\beta)$ and
$r(\Lambda_u) = 1$.

\definition{\Itemx Definition}
The system $(r_e)$ is \Def{strictly contracting} iff $r_e < 1$
for all $e$.
\enddefinition

\head 2. Perron--Frobenius Theory \endhead\newsection

The computations will use some information from the Perron--Frobenius
theory of nonnegative matrices.  This material may be found
in \cite{\gantmacher}, \cite{\minc}, or \cite{\seneta}.
We will state here the information that is used below.

Let $A$ be an $n \times n$ matrix.  We say that $A$ is nonnegative,
and write $A \ge 0$, if all entries of $A$ are nonnegative.
When $A$ and $B$ are both $n \times n$ matrices, we write
$A \ge B$ if $A - B \ge 0$.  Similarly, if $\boldx$ is an\
$n \times 1$ column vector we write $\bold x \ge 0$ if
all entries are nonnegative, and we write $\bold x > 0$
of all entries are positive.

The \Def{spectral radius} of an $n \times n$ matrix $A$ is the maximum
absolute value of the complex eigenvalues of $A$.  Write $\rho(A)$
for the spectral radius of $A$.

An $n \times n$ matrix $A$ is \Def{reducible} iff the index set
$\{1, 2, \cdots, n\}$ can be partitioned into two nonempty sets
$I, J$ such that for all $i \in I$, $j \in J$,
the entry of $A$ in row $i$ column $j$ is zero.
The matrix $A$ is \Def{irreducible} iff it is not reducible.

We cite here a few results from the Perron--Frobenius theory of
nonnegative matrices:

\Itemx If $A \ge 0$ is irreducible, then the spectral radius
$\lambda = \rho(A)$ is an eigenvalue of $A$ and there
is a positive eigenvector $\boldx > 0$ with $A\boldx = \lambda \boldx$.
\cite{\gantmacher\rm, Theorem 2, p.~53}, \cite{\minc\rm, Theorem 4.1, p.~11},
\cite{\seneta\rm, Theorem 1.5, p.~20}.

\Itemx If $A \ge 0$, $A$ irreducible, $\boldx \ge 0$, $\boldx \ne 0$,
$A\boldx = \lambda\boldx$, then $\lambda$ is the spectral
radius $\rho(A)$.
\cite{\gantmacher\rm, Remark 3, p.~63}, \cite{\minc\rm, Theorem 4.4, p.~16},
\cite{\seneta\rm, Theorem 1.6, p.~20}.

\Itemx If $A \ge B \ge 0$, where $A$ is irreducible,
then $\rho(A) \ge \rho(B)$.  \cite{\gantmacher\rm, Lemma 2, p.~57},
\cite{\minc\rm, Corollary 2.2, p.~38}, \cite{\seneta\rm, Exercise 1.9, p.~25}.

\subhead Dimension of a Graph\endsubhead
The ``dimension'' for a Mauldin--Williams graph
is computed as follows (see \cite{\mw}, \cite{\edgar}).

\definition{\Itemx Definition}
Let $\big(V, E, (r_e)_{e \in E}\big)$ be a strictly contracting
Mauldin--Williams graph.
For each positive real number $s$, let $M(s)$ be the matrix
(with rows and columns indexed by $V$) with entry
$$
	M_{uv}(s) = \sum_{e \in E_{uv}} r_e^s
$$
in row $u$ column $v$.  Let $\Phi(s) = \rho\big(M(s)\big)$
be the spectral radius
of $M(s)$.  Then $\Phi$ is continuous, strictly decreasing,
$\Phi(0) \ge 1$, $\lim_{s \to \infty}\Phi(s) = 0$.  The unique
solution $s_1 \ge 0$ of $\Phi(s_1) = 1$ is the \Def{dimension}
associated to the Mauldin--Williams graph $\big(V, E, (r_e)_{e \in E}\big)$.
\enddefinition

Note that, as a consequence of the rule for matrix
multiplication, if $k \in \N$, then the $k$th power $M(s)^k$ has entry
$$
	\sum_{\alpha \in E_{uv}^{(k)}} r(\alpha)^s
$$
in row $u$ column $v$.  Because of our conventions about
empty strings, this is true even for $k=0$.

Suppose that $(V, E)$ is a strongly connected graph.  Since
$M(s_1)$ has non-zero entry in row $u$ column $v$ whenever
there is an edge in $E_{uv}$, it follows that the matrix
is irreducible.
But $1$ is an eigenvalue, so by (\?{2.1}) there is a positive eigenvector:
$$
	\lambda_u > 0,\qquad
	\lambda_u = \sum_{v \in V} \sum_{e \in E_{uv}} r_e^{s_1}\;\lambda_v,
	\qquad \sum_{u \in V}\lambda_u = 1 .
\Itemd
$$
(For convenience, we will write $\lambdamin = \min_v \lambda_v$
and $\lambdamax = \max_v \lambda_v$. 
Thus $0 < \lambdamin \le \lambdamax \le 1$.)

This summation property shows that if we define
$$
	\mu_u\big([\alpha]\big) = r(\alpha)^{s_1}\lambda_v\qquad
	\text{for all $\alpha \in E_{uv}^{(*)}$,}
\Itemd
$$
then we have the consistency property
$$
	\mu_u\big([\alpha]\big) = \sum_{e \in E_v} \mu_u\big([\alpha e]\big)
$$
when $\alpha \in E_{uv}^{(*)}$.  Thus the $\mu_u$ extend to Borel measures
on the spaces $E_u^{(\omega)}$.

\head 3. The Proofs\endhead\newsection

Let $(V,E)$ be a strongly connected directed multigraph,
let $(S_v)_{v \in V}$ be a family
of nonempty complete metric spaces, and let $(f_e)_{e \in E}$ be a
family of functions $f_e \takes S_v \to S_u$ if $e \in E_{uv}$.
Suppose each of these functions $f_e$ satisfies a Lipschitz condition
$$
	d\big(f_e(x),f_e(y)\big) \le r_e d(x,y)
\Itemd
$$
for $x,y \in S_v, e \in E_{uv}$.  We assume all $r_e < 1$
(although \cite{\edgar\rm, Exercise 4.3.9} it is actually
enough to assume that the system
of ratios $(r_e)$ is ``contracting'' in the sense that all cycles have
ratio $< 1$).   It follows that there is a unique list
of nonempty compact sets $K_v \subseteq S_v$ satisfying
the invariance condition
$$
	K_u = \bigcup\Sb v \in V\\ e \in E_{uv}\endSb f_e\big[K_v\big]
$$
for all $u \in V$.
For future use, we define
$$
	\dmin = \min \SET{\diam K_u}{u \in V},\qquad
	\dmax = \max \SET{\diam K_u}{u \in V} .
$$

Following Definition \?{2.4}, let $s_1$ be the dimension
associated with the Mauldin--Williams
graph $\big(V, E, (r_e)_{e \in E}\big)$.  Write $M(s)$ and
$\Phi(s)$ as in that definition.

We will show (Theorem \?{3.5}) that the upper box dimensions satisfy
$$
	\overline{\dim}_{\roman B}\, K_v \le s_1
$$
for all $v$.  (See, for example, \cite{\falconer\rm, \S 9.3} for
a proof of this in the non-graph case and in Euclidean space.)

Now suppose there are also lower bounds of the form
$$
	d\big(f_e(x),f_e(y)\big) \ge r'_e \,d(x,y)
\Itemd
$$
for $x,y \in S_v, e \in E_{uv}$.
Let $s_2$ be the dimension associated with the Mauldin--Williams
graph $\big(V, E, (r'_e)_{e \in E}\big)$.
Write $M'(s)$ for the matrix and $\Phi'(s)$ for
the spectral radius in Definition \?{2.4}.

Under the assumption of the ``strong open set condition'',
we will show (Theorem \?{3.14}) that the Hausdorff dimensions satisfy
$$
	\dim_{\roman H} K_v \ge s_2
$$
for all $v$.  (See, for example, \cite{\falconer\rm, \S 9.3} for
a proof of this in the non-graph case, in Euclidean space,
and in the disjoint case.)

Certainly $r'_e \le r_e$ (since each $S_v$ has at least two points,
as we have assumed).  So $M'(s) \le M(s)$ and by (\?{2.3})
$\Phi'(s) \le \Phi(s)$ so that $s_2 \le s_1$.
The two estimates together will give us
$$
	s_2 \le \dim K_v \le s_1
$$
where ``$\dim$'' is any of the fractal dimensions mentioned above.

\subhead Upper Bounds\endsubhead
For the proof of the upper bound, we need to estimate the counts
$N_\delta(K_u)$.  We will use ``cross-cuts'' of our forest of finite
strings for this.

\definition{\Itemx Definition}
A \Def{cross-cut} is a finite set $T \subseteq E^{(*)}$
such that, for every infinite string $\sigma \in E^{(\omega)}$
there is exactly one $n$ with the restriction $\sigma\rest n \in T$.
\enddefinition

If $T$ is a cross-cut, we will write $T_u = T \cap E_u^{(*)}$
and $T_{uv} = T\cap E_{uv}^{(*)}$ for $u, v \in V$.  Then, according to
the definition,
for each $u \in V$, the cylinders $\SET{[\alpha]}{\alpha \in T_u}$
constitute a partition of $E_u^{(\omega)}$.

Let $\mu_u$ be the measures on the spaces $E_u^{(*)}$ as
in (\?{2.6}).  If $T$ is a cross-cut we have
$$
	\lambda_u = \mu_u\big(E_u^{(\omega)}\big)
	= \sum_{\alpha \in T_u} \mu\big([\alpha]\big)
	= \sum_{v \in V} \sum_{\alpha \in T_{uv}} r(\alpha)^{s_1} \lambda_v .
\Itemd
$$

\proclaim{\Itemx Theorem} Let $f_e$ satisfy upper bounds \rom{(\?{3.1}).}
Let $s_1$ be the dimension for the Mauldin--Williams graph
$\big(V, E, (r_e)_{e \in E}\big)$.
Then the attractors $K_v$ satisfy
$\dimB K_u \le s_1$.
\endproclaim 
\demo{Proof}
Let $\delta > 0$ be fixed.  Then
$$
	T = \SET{\alpha}{u,v \in V, \alpha \in E_{uv}^{(*)},
	  r(\alpha) < \delta \le r(\alpha^{-})}
\Itemd
$$
is a cross-cut (called ``first time less than $\delta$'').
For each $\alpha \in T$, we have
$\delta \rmin \le r(\alpha) < \delta$.  We may estimate
the cardinalities of the sets $T_{u}$:
$$
	\lambda_u 
	= \sum_{v \in V} \sum_{\alpha \in T_{uv}} r(\alpha)^{s_1} \lambda_v
	\ge (\delta \rmin)^{s_1} \lambdamin \,\# T_{u}
$$
and therefore, for every $u \in V$, we have
$$
	\# T_{u} \le (\delta \rmin)^{-s_1}
	\frac{\lambdamax}{\lambdamin} .
$$

But of course the set $K_u$ is covered as follows:
$$
	K_u \subseteq \bigcup_{\alpha \in T_{u}}\,\alpha K_v,
$$
and the covering sets have diameters
$$
	\diam \alpha K_v \le r(\alpha) \diam K_v
	\le \delta \dmax.
$$
Therefore we have
$$
	N_{\delta \dmax}(K_u) \le (\delta \rmin)^{-s_1}
	\frac{\lambdamax}{\lambdamin} 
$$
or (write $\eta = \delta \dmax$)
$$
	N_{\eta}(K_u) \le \left(\frac{\eta \rmin}{\dmax}\right)^{-s_1}
	\frac{\lambdamax}{\lambdamin} .
$$
Taking logarithms, dividing by the positive number $- \log \eta$,
and letting $\eta \to 0$, we get
$$
	\dimB K_u = \limsup_{\eta \to 0} \frac{\log N_{\eta}(K_u)}{-\log\eta}
	\le s_1 .\qquad\qed
$$
\enddemo

\subhead Lower Bound: Disjoint Case\endsubhead
Now we will consider the lower bound estimate on the dimension
of the sets $K_v$ that we get from the inequalities
(\?{3.2}).  Note that we continue to
assume upper bounds (\?{3.1}) as before, so that the $f_e$ are continuous, the
images $\alpha K_u$ are compact, and so on.
(Thus $r'_e < 1$ for all $e$.)

Now as before,
$1$ is an eigenvalue of the matrix $M'(s_2)$,
so by (\?{2.1}) there is an eigenvector with positive entries:
$$
	\lambda'_u > 0,\qquad
	\lambda'_u = \sum_{v \in V} \sum_{e \in E_{uv}} {r'_e}^{s_1}\;\lambda'_v,
	\qquad \sum_{u \in V}\lambda'_u = 1 .
\Itemd
$$
(As usual, we will write $\lambdaminp = \min_v \lambda'_v$
and $\lambdamaxp = \max_v \lambda'_v$.)
If $\mu'_u$ are defined by
$\mu'_u\big([\alpha]\big) = r'(\alpha)^{s_2}\lambda'_v$
for all $\alpha \in E_{uv}^{(*)}$, then they
extend to Borel measures on $E_u^{(\omega)}$.

\definition{\Itemx Definition}We will say that the graph-directed
iterated function system $(f_e)$ with attracting sets $(K_v)$
falls in the \Def{disjoint case} if, for any $u, v, v' \in V$,
$e \in E_{uv}$, $e' \in E_{uv'}$, $e \ne e'$, we have
$f_e\big[K_v\big] \cap f_{e'}\big[K_{v'}\big] = \varnothing$.
\enddefinition

\proclaim{\Itemx Theorem}
Let $f_e$ satisfy lower bounds \rom{(\?{3.2}).}
Let $s_2$ be the dimension for the Mauldin--Williams graph
$\big(V, E, (r'_e)_{e \in E}\big)$.
Suppose the IFS falls in the disjoint case.  Then
the Hausdorff dimensions satisfy $\dimH K_u \ge s_2$.
\endproclaim 
\demo{Proof}
Because of the disjointness, and because of the compactness
of the sets $K_v$, there is $\eta > 0$ so that
$\dist\big(f_e[K_v] , f_{e'}[K_{v'}]\big) > \eta$,
for all $u, v, v', e, e'$ as in the definition.
We claim that for any incomparable $\alpha, \alpha' \in E_u^{(*)}$,
$$
\dist\big(\alpha K_v, \alpha' K_{v'}\big) > r'(\alpha^{-}) \eta .
$$
To see this,
let $\gamma$ be the longest common prefix of $\alpha, \alpha'$.
Because they are incomparable, they are both strictly longer than $\gamma$.
(Thus $\gamma \le \alpha^-$, so $r'(\gamma) \ge r'(\alpha^-)$.)
There are distinct letters $e, e'$ so that $\alpha = \gamma e \beta$,
$\alpha' = \gamma e' \beta'$.  Say $\beta \in E_{wv}^{(*)}$,
$\beta' \in E_{w'v'}^{(*)}$; these strings may be empty.
Then $\beta K_v \subseteq K_w$ and $\beta' K_{v'}\subseteq K_{w'}$;
since $e \ne e'$ we have
$\dist(e \beta K_v, e' \beta' K_{v'}) \ge \dist(e K_w, e' K_{w'}) > \eta$;
so finally $\dist(\gamma e \beta K_v, \gamma e' \beta' K_{v'})
\ge r'(\gamma) \eta \ge r'(\alpha^{-}) \eta$.

Let $B \subseteq S_u$ be a Borel set.  Let 
$\delta = \diam B / \eta$
and let $T$ be the cross-cut ``first time less than $\delta$''
$$
	T = \SET{\alpha}{u \in V, v \in V, \alpha \in E_{uv}^{(*)},
	  r'(\alpha) < \delta \le r'(\alpha^{-})} .
\Itemd
$$
The elements $\alpha \in T_u$ are all incomparable, and
$\diam B = \eta \delta \le \eta r'(\alpha^{-})$, so
at most one set $\alpha K_v$ with $\alpha \in T_u$ can meet $B$.  This means
$h^{-1}[B] \subseteq [\alpha]$ for some $\alpha \in T$,
and thus $\mu'_u\big(h^{-1}[B]\big) \le \mu'_u\big([\alpha]\big)
= r'(\alpha)^{s_2} \lambda'_v \le \lambdamaxp \eta^{-s_2}(\diam B)^{s_2}$.

Now if $\scr C$ is a countable cover of $K_u$ by
Borel sets, we have
$$
	\sum_{B \in \scr C} \left(\diam B\right)^{s_2}
	\ge \left(\frac{\eta^{s_2}}{\lambdamaxp}\right)
	\sum \mu'_u\big(h^{-1}[B]\big)
	\ge \left(\frac{\eta^{s_2}}{\lambdamaxp}\right)\mu'_u\big(K_u\big).
$$
It follows that $\scrH^{s_2}(K_u) \ge
\left({\eta^{s_2}}/{\lambdamaxp}\right)\mu'_u\big(K_u\big) > 0$.
Therefore $\dimH K_u \ge s_2$.
\qed\enddemo

\subhead Lower Bound: Strong Open Set Condition\endsubhead
To prove a lower bound for the dimension
of an IFS we must of course limit the overlap.  This is often done with
an ``open set condition''.  In general complete metric spaces
(not necessarily subsets of Euclidean space) we require
the ``strong open set condition''.   For the strong open
set condition, and its application in this calculation,
we follow \cite{\bandtg}, \cite{\patzschke}, \cite{\schief}.

\definition{\Itemx Definition}  Let $(V,E)$ be a directed multigraph,
let $(S_v)_{v \in V}$ be complete metric spaces, let
$f_e \takes S_v \to S_u$ for $e \in E_{uv}$ be a strictly contracting
graph-directed IFS, and let $(K_v)_{v \in V}$ be the attractor
list.  We say that the IFS $(f_e)$ satisfies the \Def{open set condition}
iff, for each $v \in V$ there is an open set $U_v \subseteq S_v$
such that:
\roster
\item"(a)" For all $u, v \in V$, $e \in E_{uv}$,
  we have $f_e\big[U_v\big] \subseteq U_u$.
\item"(b)" For all $u, v, v' \in V$, $e \in E_{uv}$, $e' \in E_{uv'}$,
  $e \ne e'$,
  we have $f_e\big[U_v\big] \cap f_{e'}\big[U_{v'}\big] = \varnothing$.
\endroster
We say that the IFS $(f_e)$ satisfies the \Def{strong open set condition}
iff there exist open sets $U_v$ satisfying (a), (b), and
\roster
\item"(c)"  For all $v \in V$, we have $U_v \cap K_v \ne \varnothing$.
\endroster
\enddefinition

Unlike the case where each $S_v$ is Euclidean space, it now no longer
follows that the images $f_e\big[U_v\big]$ are open sets.  Even if
$f_e$ is a homeomorphism of $U_v$ onto $f_e\big[U_v\big]$, it need
not follow that $f_e\big[U_v\big]$ is open.

Notice the following consequence of the open set condition.
(Of course, (\?{3.2}) implies that $f_e$ is injective.)

\proclaim{\Itemx Lemma}  Let $f_e$ satisfy the open set condition
with open sets $U_v$.  Suppose the maps $f_e$ are injective.
If $\alpha \in E_{uv}^{(*)}, \beta \in E_{uv'}^{(*)}$
are incomparable, then
$f_\alpha\big[U_v\big] \cap f_{\alpha'}\big[U_{v'}\big] = \varnothing$.
\endproclaim 
\demo{Proof}
Let $\gamma$ be the longest common prefix of $\alpha, \alpha'$.
Because they are incomparable, they are both strictly longer than $\gamma$.
There are distinct letters $e, e'$ so that $\alpha = \gamma e \beta$,
$\alpha' = \gamma e' \beta'$.  Say $\beta \in E_{wv}^{(*)}$,
$\beta' \in E_{w'v'}^{(*)}$; these strings may be empty.  Now
by (a), we have
$\beta U_v \subset U_w$ and $\beta' U_{v'} \subseteq U_{w'}$.
Next, $e\beta U_v \subseteq f_e[U_w]$ and
$e'\beta' U_{v'}\subseteq f_{e'}[U_{w'}]$ and by (b) these are
disjoint.  Since the maps $f_e$ are injective, so are their compositions,
so we conclude that $\alpha U_v = \gamma e\beta U_v$ and
$\alpha' U_{v'} = \gamma e'\beta' U_{v'}$ are also disjoint.
\qed\enddemo

From the uniqueness properties of the invariant sets $K_v$,
it is easy to see that $K_v \subseteq \overline{U_v}$.  But
$K_v \subseteq U_v$ need not be true.
The strong open set condition provides a substitute for this.

\proclaim{\Itemx Lemma} Let $(V,E)$ be a strongly connected directed multigraph.
Suppose the IFS $(f_e)$ satisfies the strong open set condition
with open sets $U_v$.
For each $v \in V$, there is a cycle $\zeta \in E_{vv}^{(*)}$
such that $\zeta K_v \subseteq U_v$.
\endproclaim 
\demo{Proof}
By (c), there is a point $x \in K_v \cap U_v$.  Let $\sigma \in E_v^{(\omega)}$
be an address of $x$.  Then the singleton
$\{x\}$ is the decreasing intersection of compact sets
of the form $(\sigma\rest n) K_{u_n}$.  Now $U_v$ is a neighborhood
of the point $x$, so we conclude by compactness that there is
$n$ so that $(\sigma\rest n) K_{u_n} \subseteq U_v$.  Because the
graph $(V,E)$ is strongly connected, there is a path
$\gamma \in E_{u_nv}^{(*)}$.  Let $\zeta = (\sigma\rest n) \gamma$.
Then we have $\gamma K_v \subseteq K_{u_n}$ so that
$\zeta K_v = (\sigma\rest n) \gamma K_v \subseteq
(\sigma\rest n)K_{u_n} \subseteq U_v$.
\qed\enddemo

\proclaim{\Itemx Theorem} Let $f_e$ satisfy lower bounds \rom{(\?{3.2}).}
Let $s_2$ be the dimension for the Mauldin--Williams graph
$\big(V, E, (r'_e)_{e \in E}\big)$.
Suppose $(f_e)$ satisfies the strong open
set condition.  Then the attractors $K_v$ satisfy
$\dimH K_v \ge s_2$.
\endproclaim 
\demo{Proof}
Let $s < s_2$.  We will show that $\dimH K_v \ge s$.  Now
$\Phi'(s) > 1$, so the spectral radius $t = \rho\big(M'(s)\big)$ is $> 1$.

Let $U_v$ be open sets for the strong open set condition.
For each $v \in V$, choose a cycle $\zeta_v \in E_{vv}^{(*)}$
so that $\zeta_v K_v \subseteq U_v$.  Then
$$
	c = \min_{v \in V} r'(\zeta_v)^s
$$
is a positive number.  Let $n \in \N$ be such that $c t^n > 1$.

Now we will define a new IFS.  First we define a
new directed multigraph $(V, \widetilde E)$.
The vertices of the new graph will be $V$ as before.
The edges from vertex $u$ to vertex $v$ will be
$$
	\widetilde E_{uv} = \SET{\alpha \zeta_v}{\alpha \in E_{uv}^{(n)}} .
$$
The corresponding maps are $\tilde f_{\alpha \zeta_v} = f_{\alpha \zeta_v}$.
Note that $\alpha \zeta_v [K_v] \subseteq K_u$, so the new attractors
$\widetilde K_v$ are subsets of the old attractors $K_v$.

We claim that that the new IFS falls in the disjoint case.  Indeed,
let $\alpha \zeta_v \in \widetilde E_{uv}$ and
$\alpha' \zeta_{v'} \in \widetilde E_{uv'}$ be distinct edges.
Then $\zeta_v K_v \subseteq U_v$ and $\zeta_{v'} K_{v'} \subseteq U_{v'}$.
But $\alpha$ and $\alpha'$ have the same length and are unequal, so they
are not comparable, and thus by Lemma \?{3.12} we have
$\alpha \zeta_v K_v \cap \alpha' \zeta_{v'} K_v'
\subseteq \alpha U_v \cap \alpha' U_{v'} = \varnothing$.

Now the matrix $\widetilde M'(s)$ for the new IFS has entry
$$
	\sum_{\alpha \in E_{uv}^{(n)}} r'(\alpha \zeta_v)^s
$$
in row $u$ column $v$.  But this entry is 
$$
	\ge\, c \sum_{\alpha \in E_{uv}^{(n)}} r'(\alpha)^s
$$
which is the corresponding entry in the matrix $c M'(s)^n$.
Now by (\?{2.3}),
$$
	\rho\big(\widetilde M'(s)\big)
	\ge \rho\big(c M'(s)^n\big) = c t^n > 1 .
$$
This shows (by Theorem \?{3.9}) that $\dimH \widetilde K_v > s$,
so also $\dimH K_v > s$.
\qed\enddemo

\head 4. Example: A Julia Set \endhead\newsection

We will consider an example showing how the theorem may be
applied to estimate the fractal dimension of a set.  This example
is the Julia set for the transformation $z^2 - 1/2$
of the complex plane $\Bbb C$.  See, for example,
\cite{\falconer\rm, Chapter~14}
for an explanation of Julia sets in simple cases like this.
McMullen \cite{\mcmullen} has provided an algorithm that
allows computation of Hausdorff dimension up to any desired
accuracy of conformal expanding open maps (so having Markov partition) 
which include this example.

The Julia set we will consider here is shown in Figure \?{4.1}.
See \cite{\golds} for more details of the following computation.

  \global\advance\itemno by 1%
  \topinsert
\centerline{\BoxedEPSF{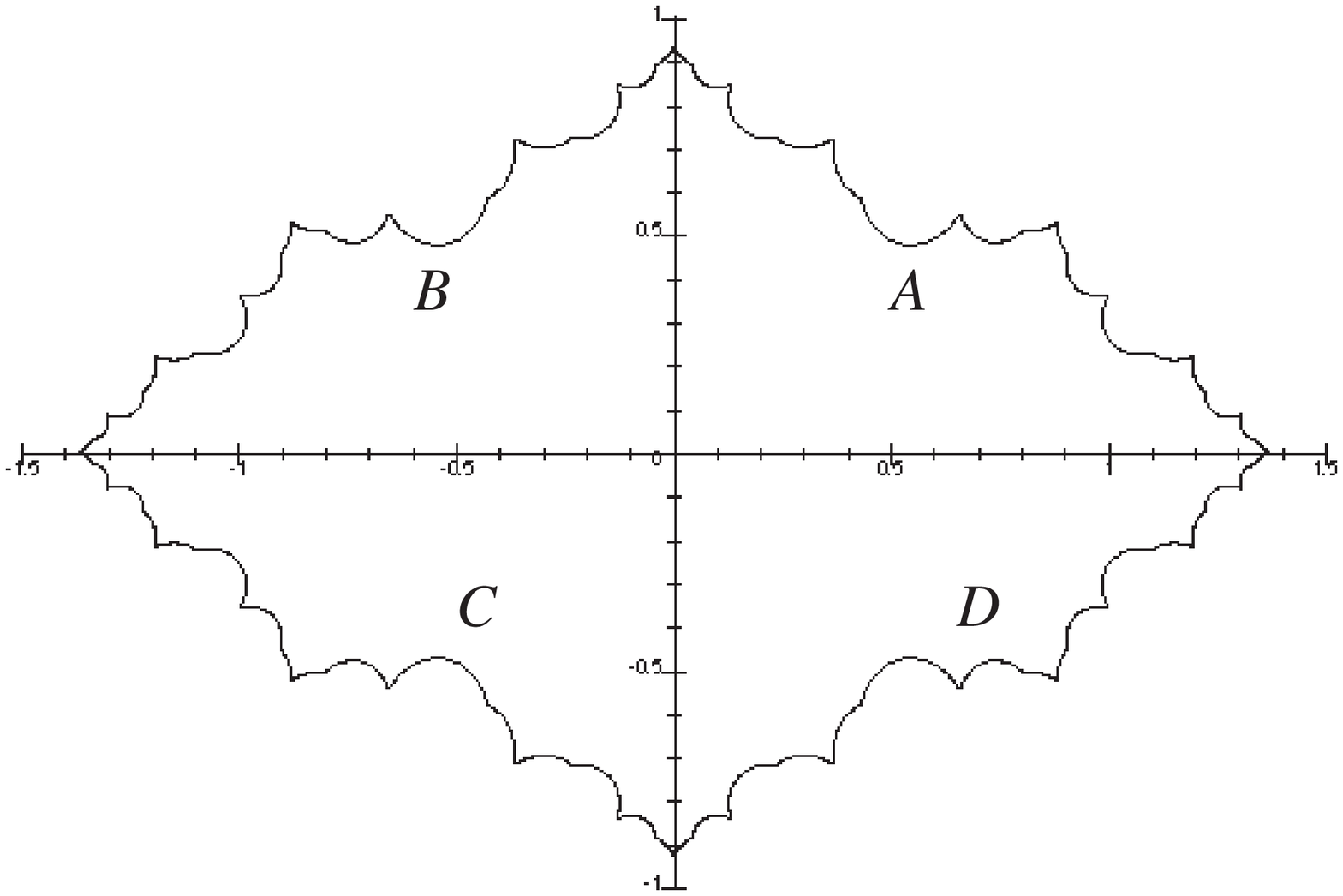 scaled 400}}
\botcaption{Figure \theitem}A Julia set\endcaption
  \endinsert

We divide the set into four parts $A, B, C, D$, using the four quadrants
of the plane.  Then under the two inverse maps $\pm\sqrt{w+1/2\,}$
each of these parts is made up of the images of two of the parts.
$A$ is made up of images of $A$ and $B$;
$B$ is made up of images of $C$ and $D$;
$C$ is made up of images of $A$ and $B$;
$D$ is made up of images of $C$ and $D$.
Thus our directed graph will have four vertices and eight
edges.  See Figure \?{4.2}.

  \global\advance\itemno by 1%
  \topinsert
\centerline{\BoxedEPSF{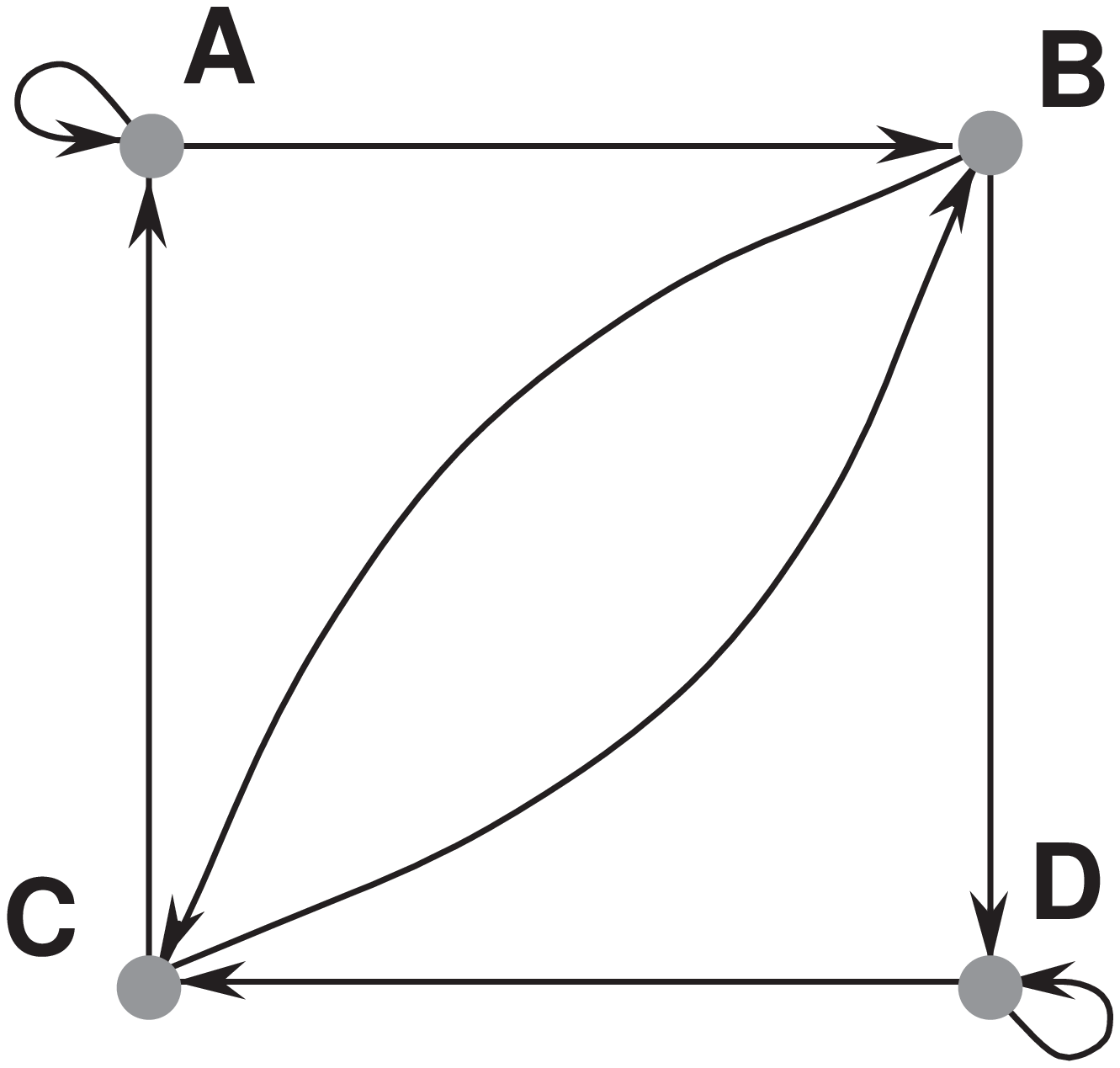 scaled 250}}
\botcaption{Figure \theitem}The graph\endcaption
  \endinsert

For the upper and lower bounds, we will use a simple estimate
based on the derivatives.

\proclaim{\Itemx Lemma}
Suppose $\varphi(z) = z^2 - 1/2$ maps convex set $U \subseteq \Bbb C$
bijectively onto
convex set $W$, and $0\not\in W$.  Write $f(w) = \sqrt{w+1/2\,}$
for the branch of the inverse that maps $W$ onto $U$.  Let
$$
	m = \inf\SET{|z|}{z \in U};\qquad M = \sup\SET{|z|}{z \in U}.
$$
Then, for any $w_1, w_2 \in W$,
$$
	\frac{1}{2 M}\;\big|w_2-w_1\big|
	\le \big|f(w_2)-f(w_1)\big|
	\le \frac{1}{2 m}\;\big|w_2-w_1\big| .
$$
\endproclaim 
\demo{Proof}  Note that $\varphi'(z) = 2 z$, so $|\varphi'(z)| \le 2 M$
for all $z \in U$.  Also $f'(w) = 1/(2 f(w))$ and if
$w \in W$ then $z = f(w) \in U$, so $|f'(w)| \le 1/(2m)$ for
all $w \in W$.

Now let $w_1, w_2 \in W$, and write $z_1 = f(w_1)$,
$z_2 = f(w_2)$.  Since $U$ is convex
we may integrate along the line segment joining $z_1$ to $z_2$:
$$
\align
	\big|w_2-w_1\big|
	&= \big|\varphi(z_2)-\varphi(z_1)\big|
	= \left|\int_{z_1}^{z_2} \varphi'(z)\,dz\right|\\
	&\le \left(\sup_{z \in U} |\varphi'(z)|\right) \big|z_2-z_1\big|
	= 2 M \big|z_2-z_1\big| .
\endalign
$$
Consequently we have $|f(w_2)-f(w_1)| \ge \big(1/(2M)\big) |w_2-w_1|$.
Similarly,
$$
\align
	\big|f(w_2)-f(w_1)\big|
	&= \left|\int_{w_1}^{w_2} f'(w)\,dw\right|\\
	&\le \left(\sup_{w \in W} |f'(w)|\right) \big|w_2-w_1\big|
	= \frac{1}{2m}\, \big|w_2-w_1\big| .\qquad\qed
\endalign
$$
\enddemo

For the convex sets required here, we begin with the smallest circle
that surrounds our set, a certain parallelogram inside, and segments
on the imaginary axis.  We get four three-sided regions, each containing
the portion of the curve in one quadrant.  See Figure \?{4.4}.

  \global\advance\itemno by 1%
  \topinsert
\centerline{\BoxedEPSF{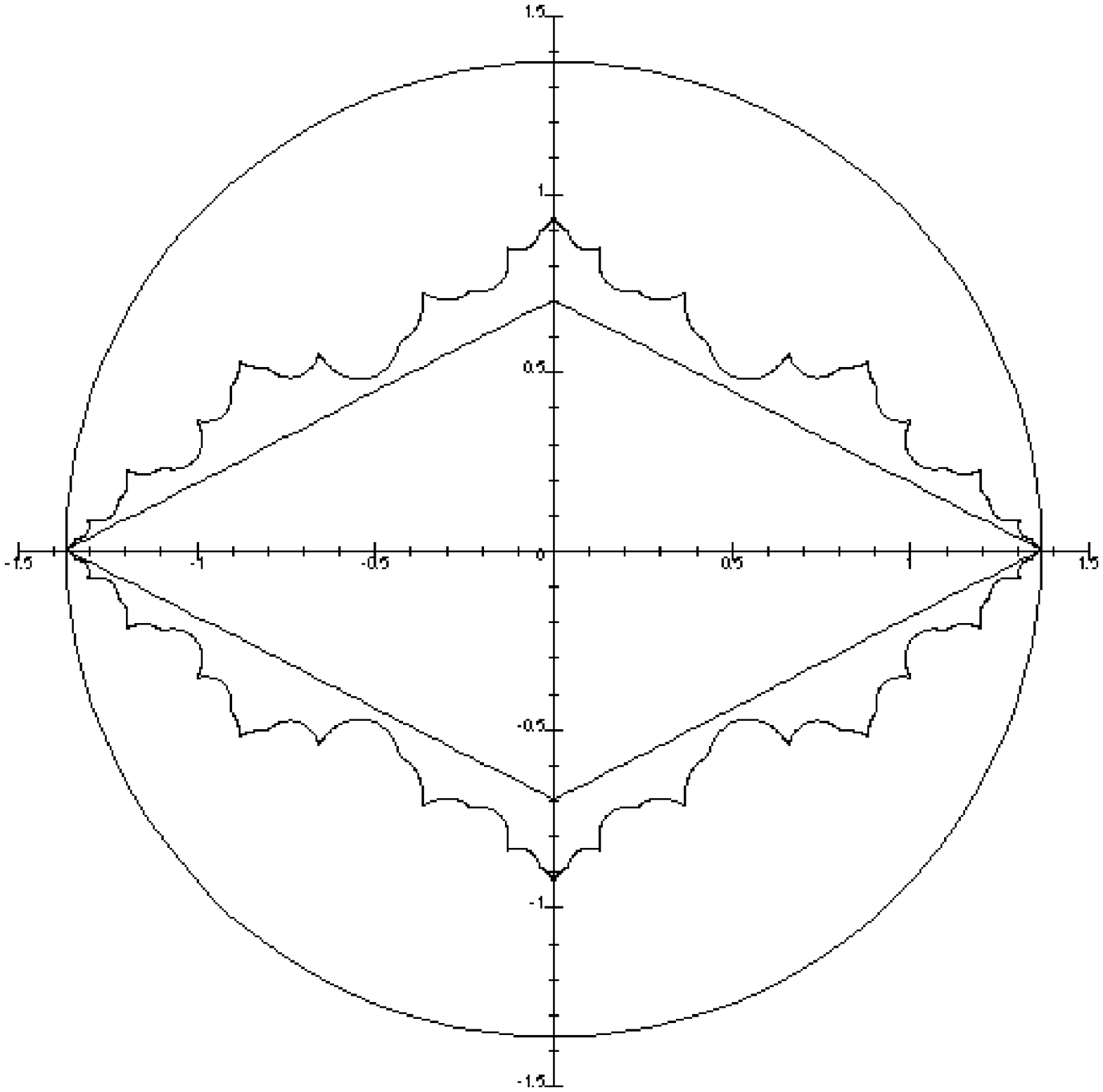 scaled 400}}
\botcaption{Figure \theitem}Inner and outer bounds\endcaption
  \endinsert

The parallelogram is chosen so that for each of the four regions bounded
by our curves, the images under
$\pm\sqrt{w+1/2\,}$ are inside the original regions.
(See Figure \?{4.5}.)

  \global\advance\itemno by 1%
  \topinsert
\centerline{\BoxedEPSF{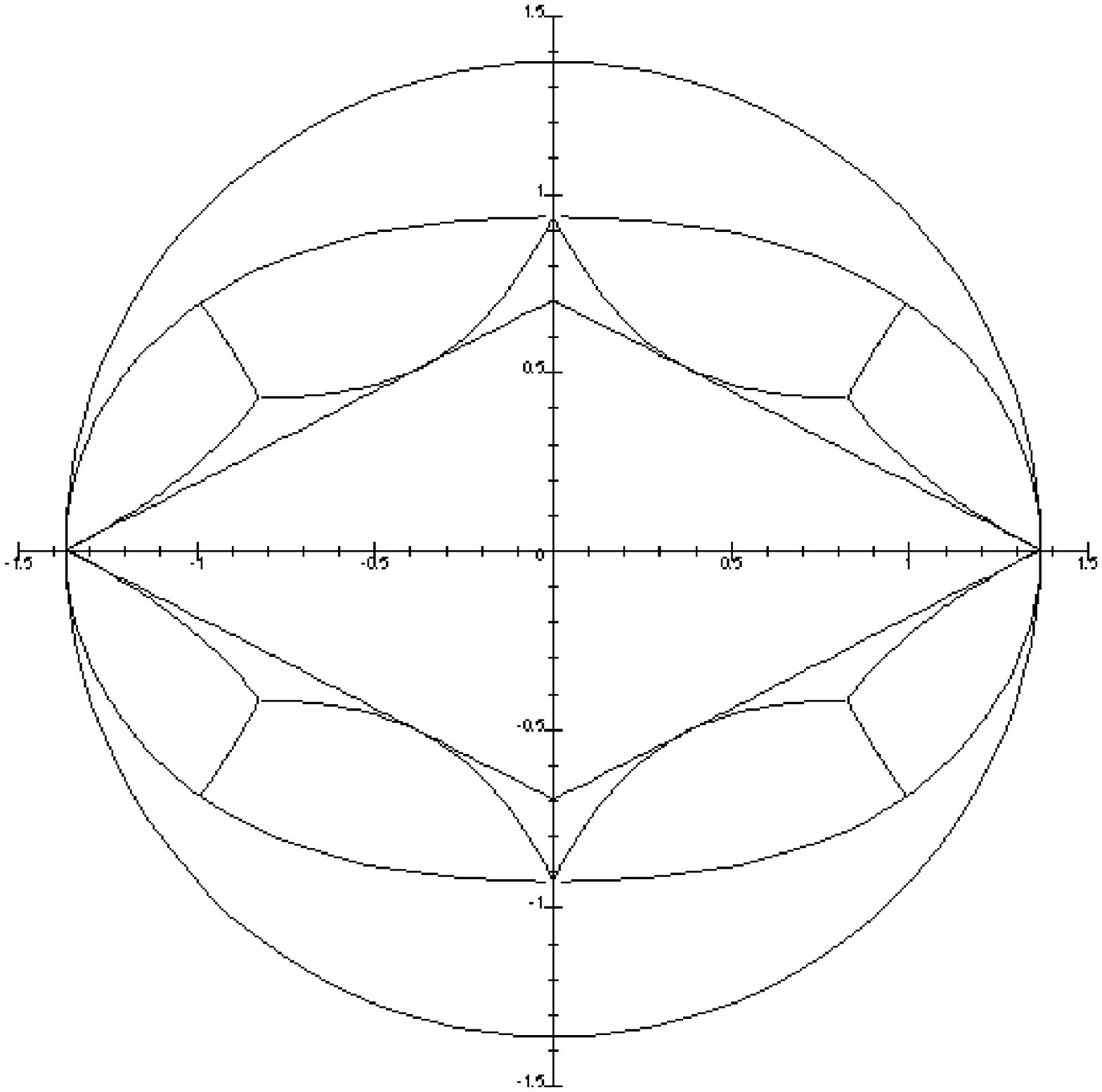 scaled 400}}
\botcaption{Figure \theitem}Images are inside the regions\endcaption
  \endinsert

Now we are in a position to estimate the constants for the
inequalities.  The outer circle has radius $(1 + \sqrt{3}\,)/2$,
so $M = (1 + \sqrt{3}\,)/2$ and our lower bounds are all
$r'(e) = 1/(1+\sqrt{3}\,)$.  Then
$$
	M'(s) = \bmatrix
	a & a & 0 & 0\\
	0 & 0 & a & a\\
	a & a & 0 & 0\\
	0 & 0 & a & a\\
	\endbmatrix
$$
where $a = (1+\sqrt{3})^{-s}$.  The spectral radius is $1$
when $s = s_2 = \log 2/\log(1+\sqrt{3}\,) > 0.689$.

Similarly, for the closest point we get
$m = \sqrt {9+3\,\sqrt {3}}/6$, and $r(e) = 1/(2m) \approx 0.7962$.
Computing the spectral radius as before, we get $s_1 < 3.042$.

The (rather poor) estimates are
$$
	0.689 < \dim K < 3.042\,.
$$

But that need not be the end of our work.  Consider the eight small
regions in Figure \?{4.5}.  Their images, in turn, fall as shown in
Figure \?{4.6}.  This graph-directed iterated function system has
$8$ nodes and $16$ edges.  We may work with an $8 \times 8$ matrix.
[Actually, because of symmetry, we may reduce this to a $2 \times 2$
matrix, but such questions of efficiency are not dealt with in
this paper.]  Using the upper and lower estimates for these regions,
we obtain
$$
	0.735 < \dim K < 1.758 .
$$
A considerable improvement, but still not very good.

  \global\advance\itemno by 1%
  \topinsert
\centerline{\BoxedEPSF{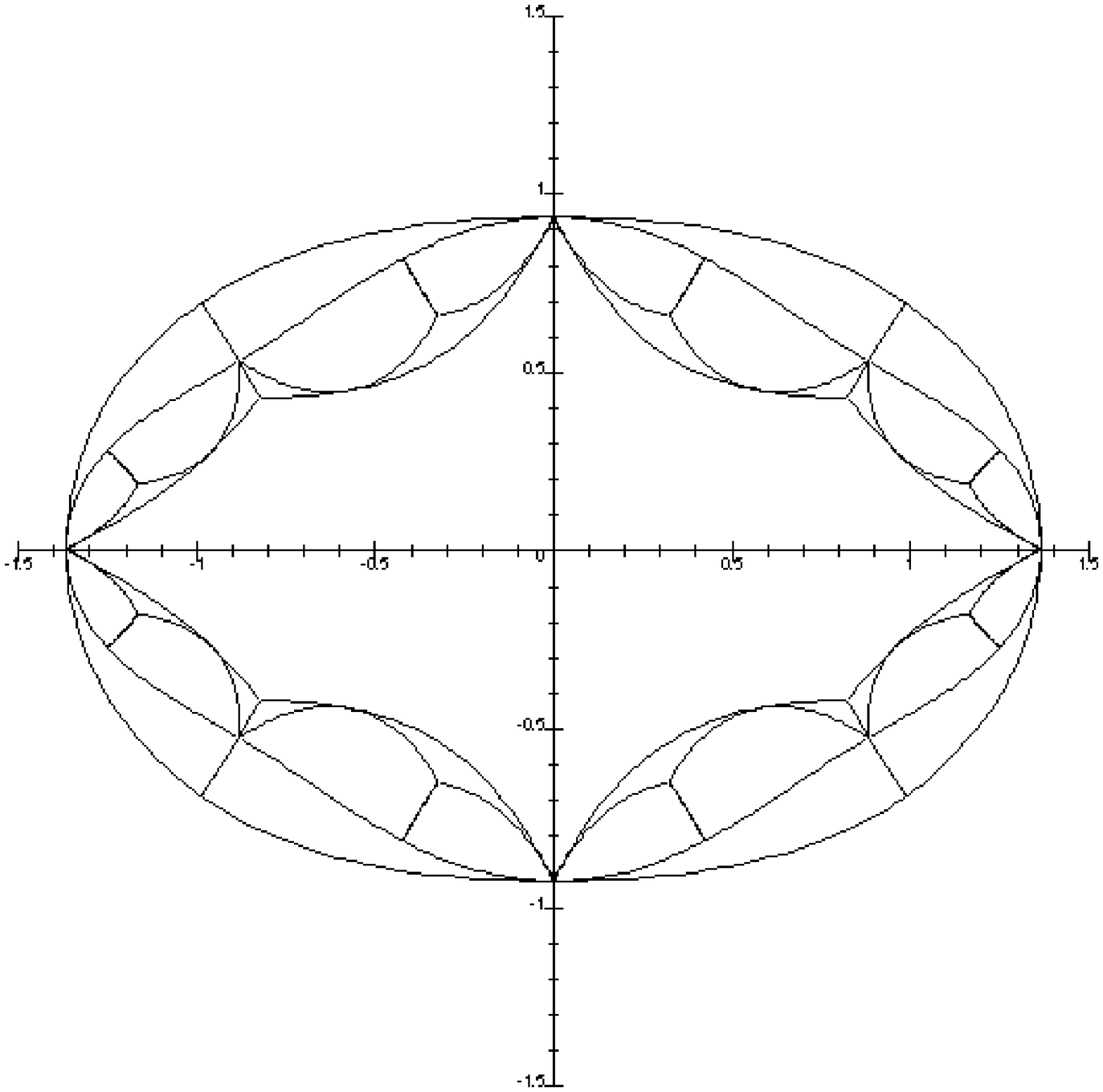 scaled 400}}
\botcaption{Figure \theitem}Eight regions and their images\endcaption
  \endinsert

If we have a computer to help us, we may continue to get better
approximations.  When we reach a matrix of size $2^{12} \times 2^{12}$
our estimate is
$$
	1.069 < \dim K < 1.077\;.
$$
The size of the computation is not as bad as the size of the
matrix may suggest: each row has only two nonzero entries.

McMullen's method \cite{\mcmullen} computes the dimension for
this Julia set as $1.07336$ (according to the computer program
on his web site).  This agreement provides evidence that our
methods (and our computer program) are correct.

Other Julia sets for the quadratic maps $z^2+c$ are considered
in \cite{\golds}.  For $c$ outside the Mandelbrot set, and for
$c$ in the main cardioid of the Mandelbrot set, the process
required to carry out the inner and outer approximations
was automated.  It is hoped that details of this will be
published elsewhere.

Note that the functions $z^2+c$ used in this example are ``conformal'' in
the sense used in geometry (except at $z=0$).
So this example is a ``self-conformal''
set.  See \cite{\mauldin} and \cite{\strichartz} for more information
on self-conformal fractals.

For a map $f \takes S \to T$ of metric spaces, one might say that
$f$ is \Def{conformal} iff, for every $x_0 \in S$, there is a positive
constant $c(x_0)$ such that for every $\varepsilon > 0$
there is $\delta > 0$ such that
$$
	\left|\frac{d\big(f(x),f(x_0)\big)}{d(x,x_0)} - c(x_0)\right| < \varepsilon
$$
whenever $d(x,x_0) < \delta$.  But as far as we know, the only cases
that have been considered (aside from similarities in general
metric spaces) are differentiable manifolds.

\enddocument